\documentclass[12pt]{article}
\usepackage[cp1251]{inputenc}
\usepackage{amssymb,amsmath,amscd}
\usepackage[all]{xy}
\usepackage[english]{babel}

\newtheorem{Th}{Theorem}

\newtheorem{Ex}{Example}

\newtheorem{Rem}{Remark}

\begin{document}

\begin{center}
\textsc{Cartesian products of the $g$-topologies are a $g$-topology}\\
\medskip
\textsl{Jumaev Davron Ilxomovich}\footnote{\scriptsize Tashkent Institute of Architecture and Civil Engineering, 7, Kichik Xalqa yuli Str.,  Tashkent 100084, Uzbekistan,\\
e-mail: d-a-v-ron@mail.ru},\\
\textsl{Ishniyazov Baxrom Normamatovich}\footnote{\scriptsize Tashkent State Agrarian University, 2, University Str., Kibray District, Tashkent Region 100700, Uzbekistan,\\
e-mail: ishniyazov73@mail.ru},\\
\textsl{Tagaymuratov Abror Olimovich} \footnote{\scriptsize Chirchik State Pedagogical Institute, 104,  Amir Temur Str., Chirchik town, Tashkent Region 111700, Uzbekistan,\\
e-mail: abror.t93@mail.ru}

\end{center}
\begin{abstract}
We show that unlike the usual topologies the $g$-topologies are closed with respect to the Cartesian products. Moreover, we bring much detailed explanations some examples of concepts related the statistical metric spaces.\\


\end{abstract}

\section{Introduction}
As mentioned in \cite{Menger1942} the three principal applications of statistical metrics are to macroscopic, microscopic and
physiological spatial measurements. Statistical metrics are designed to provide us firstly with a method removing conceptual difficulties
from microscopic physics and transferring them into the underlying geometry, secondly  with a treatment of thresholds of spatial sensation
eliminating the intrinsic paradoxes of the classical theory.

The notion of distance is defined in terms of functions, points and
sets. Indeed, in many situations, it is appropriate to look upon the distance concept as a statistical rather than a determinate one.
More precisely, instead of associating a number to the distance $d(p,q)$ with every pair of points $p$, $q$, one should associate a
distribution function $F_{pq}$ and for any positive number $x$, interpret $F_{pq}(x)$ as the probability that the distance from p to q be
less than $x$.

Using this idea, K.~Menger in \cite{Menger1942} defined a statistical metric space using the probability function in the year
1942. In 1943, shortly after the appearance of Mengers article, Wald published an article \cite{Wald1943} in which he criticized Mengers generalized triangle inequality.

In \cite{Thorp1962} the following questions raised by Thorp in statistical metric spaces:
\begin{itemize}
\item What are the necessary and sufficient conditions that the $g$-topology of type $V$ to be of type $V_{D}$?

\item What are the necessary and sufficient conditions that the $g$-topology of type $V_{\alpha}$ to be the $g$-topology of type $V_{D}$?

\item What conditions are both necessary and sufficient for the $g$-topology of type $V_{\alpha}$ to be a topology?
\end{itemize}

In \cite{RenVad2019} it had given partial answer to the above questions. Also, it was provided the basis for carrying out analysis in statistical metric spaces, in particular for the development of various $g$-topologies, neighborhoods defined in a statistical metric space and also the improvement of $\lambda_{\Omega}$-open sets in a generalized metric space. The authors of \cite{RenVad2019} had given more examples of the neighborhoods defined in a statistical metric space and the special kind of relationship between various $g$-topologies defined by Thorp in a $SM$ space. We seem that the examples have shortcomings. That is why we in the present paper completed the shortcomings and give more detail clarifies.

Further, we give positive answers to the question, expressed in \cite{Zaitov2019NUUz}, which asks: Is the Cartesian product of the $g$-topological spaces a $g$-topological space?

\section{Preliminaries}

A statistical metric space ($SM$ space) is an ordered pair $(S,\, F)$ where $S$ is a non-null set and $F$ is a mapping from $S\times S$ into the set of distribution functions (that is, real-valued functions of a real variable which are everywhere defined, non decreasing,
left-continuous and have infimum $0$ and supremum $1$).

The distribution function $F(p,\, q)$ associated with a pair of points p and q in S is denoted by $F_{pq}$. Moreover, $F_{pq}(x)$ represents the probability that the ``distance'' between $p$ and $q$ is less than $x$.

The functions $F_{pq}$ are assumed to satisfy the following:

\begin{itemize}
\item[$(SM\mbox{-}I)$] 	$F_{pq}(x)=1$ for all $x>0$ if and only if $p=q$.

\item[$(SM\mbox{-}II)$] 	$F_{pq}(0)=0$.

\item[$(SM\mbox{-}III)$] 	$F_{pq}=F_{qp}$.

\item[$(SM\mbox{-}IV)$] 	If $F_{pq}(x)=1$ and $F_{qr}(y)=1$, then $F_{pr}(x+y)=1$.
\end{itemize}

We often find it convenient to work with the tails of the distribution functions rather than with these distribution functions
themselves. Then the tail of $F_{pq}$, denoted by $G_{pq}$, is defined by $G_{pq}(x)=1-F_{pq}(x)$ for all $ x \in\mathbb{R}.$

Let $(S,\, F)$ be a statistical metric space. Then the Menger inequality is,

\begin{itemize}
\item[$(SM\mbox{-}IVm)$] 	$F_{pr}(x+y)\geq T(F_{pq}(x),\, F_{qr} (y))$
\end{itemize}
holds for all points $p$, $q$, $r \in S$ and for all numbers $x$, $y\geq 0$ where $T$ is a $2$-place function on the unit square satisfying:

\begin{itemize}
\item[$(T\mbox{-}I)$] $0 \leq T(a,\, b)\leq 1$ for all $a,\, b\in [0,1]$.

\item[$(T\mbox{-}II)$] $T(c,\, d)\geq T(a,\, b)$ for all $a,\, b,\, c,\, d \in [0,1]$ such that $c\geq a$, $d\geq b$ (monotonicity).

\item[$(T\mbox{-}III)$] $T(a,\, b)=T(b,\, a)$ for all $a, b\in [0,1]$ (commutativity).

\item[$(T\mbox{-}IV)$] $T(1,\, 1)=1$.

\item[$(T\mbox{-}V)$] $T(a,\, 1)>0$ for all $a>0$.
\end{itemize}

Let $(S,\ F)$ be a statistical metric space, $p\in S$ and $u$, $v$ be positive numbers. Then
\begin{gather*}
N_{p}(u,\, v)= \{q\in S:\, F_{pq} (u)>1-v\}=\{q\in S:\, G_{pq} (u)<v\}
\end{gather*}
is called the $(u,\, v)$-sphere with the center $p$.

The following example shows the existence of $(u,\, v)$-sphere in a statistical metric space.

\begin{Ex}\label{Npuv}
{\rm Consider the $SM$ space $(S,\, F)$ where $S$ denotes the possible outcomes of getting a tail when a coin is tossed once. Then $S=\{0,\, 1\}$. Let $F_{pq}(u)$ be the probability that the ``distance'' between $p$ and $q$ is less than u where $u>0$ and $p,\, q\in S$. We have (A.~A.~Zaitov):

$F_{00} (u)=F_{11} (u)=1$ for all $u>0$;
and
\begin{gather*}
F_{01} (u)=F_{10} (u)=
\begin{cases}
0, & \text{if}\,\, 0<u\leq 1,\\
1, & \text{if}\,\, u>1.
\end{cases}
\end{gather*}

Fix $p=0$. Then
\begin{gather*}
N_{0} (u,\, v)=
\begin{cases}
\{0\}, & \text{if}\,\, 0<u\leq 1, 0\leq v\leq 1,\\
\{0,1\}, &\text{if}\,\,   0<u\leq 1, v>1,\\
\{0,1\}, &\text{if}\,\,   u>1, v\geq 0.
\end{cases}
\end{gather*}

Fix $p=1$. Then
\begin{gather*}
N_{1} (u,\, v)=
\begin{cases}
\{1\}, & \text{if}\,\,  0<u\leq 1, 0\leq v\leq 1,\\
\{0,1\}, & \text{if}\,\,  0<u\leq 1, v>1,\\
\{0,1\}, & \text{if}\,\,   u>1, v\geq 0.
\end{cases}
\end{gather*}}
\end{Ex}

For fixed positive numbers $u$ and $v$, define a set
\begin{gather*}
U(u,\, v)=\{(p,q)\in S\times S: G_{pq} (u)<v\}.
\end{gather*}

Let us illustrate such sets in the following example.

\begin{Ex}\label{Uuv}
{\rm Consider the $SM$ space $(S,F)$ where $S$ denotes the possible outcomes of rolling a dice. Then $S=\{1,\, 2,\, 3,\, 4,\, 5,\, 6\}$ and the distribution function $F_{pq}(x)$ is the probability that the ``distance'' between $p$ and $q$ is less than u where $u>0$ and $p,\, q\in S$. Consider the usual metric on $S$ induced on the real line, i.~e. $d(p,\, q)=|q-p|$. We have (A.~A.~Zaitov):

\begin{align*}
&F_{pp}(u)=1\,\, \text{for all}\,\, u>0\,\, \text{and}\,\, p=1,\, 2,\, 3,\, 4,\, 5,\, 6;\\
&F_{p(p+1)}(u)=F_{(p+1)p}(u)=
\begin{cases}
0, &  0<u\leq 1,\\
1, &  u>1,
\end{cases} \qquad\mbox{where}\quad p=1,\, 2,\, 3,\, 4,\, 5;\\
&F_{p(p+2)}(u)=F_{(p+2)p}(u)=
\begin{cases}
0, & 0<u\leq 2,\\
1, & u>2,
\end{cases} \qquad\mbox{where}\quad p=1,\, 2,\, 3,\, 4;\\
&F_{p(p+3)}(u)=F_{(p+3)p}(u)=
\begin{cases}
0, & 0<u\leq 3,\\
1, & u>3,
\end{cases} \qquad\mbox{where}\quad p=1,\, 2,\, 3;\\
&F_{p(p+4)}(u)=F_{(p+4)p}(u)=
\begin{cases}
0, & 0<u\leq 4,\\
1, & u>4,
\end{cases} \qquad\mbox{where}\quad p=1,\, 2;\\
&F_{p(p+5)}(u)=F_{(p+5)p}(u)=
\begin{cases}
0, & 0<u\leq 5,\\
1, & u>5,
\end{cases} \qquad\mbox{where}\quad p=1.
\end{align*}

Consequently:

for every pair of $u$, $v$ such that $0<u\leq 1$, $0<v\leq 1$ we have:
\begin{gather*}
U(u,\, v)=\{(p,\, q)\in S\times S:\, d(p,\, q)=0\}=\{(p,\, p):\, p\in S\} = \Delta(S);
\end{gather*}

for every pair of $u$, $v$ such that $1<u\leq 2$, $0<v\leq 1$:
\begin{multline*}
U(u,\, v)=\{(p,\, q)\in S\times S:\, d(p,\, q)\leq 1\} = \\
=\Delta (S)\cup \{(p,\, p+1):\, p=1,\, 2,\, 3,\, 4,\, 5\}\cup \{(p+1,\, p):\, p=1,\, 2,\, 3,\, 4,\, 5\}\stackrel{def} = \Delta_{1}(S);
\end{multline*}

for every pair of $u$, $v$ such that $2<u\leq 3$, $0<v\leq 1$:
\begin{multline*}
U(u,\, v) = \{(p,\, q)\in S\times S:\, d(p,\, q)\leq 2\} =\\
=\Delta_{1}(S)\cup \{(p,\, p+2):\, p=1,\, 2,\, 3,\, 4\}\cup \{(p+2,\, p):\, p=1,\, 2,\, 3,\, 4\}\stackrel{def} = \Delta_{2}(S);
\end{multline*}

for every pair of $u$, $v$ such that $3<u\leq 4$, $0<v\leq 1$:
\begin{multline*}
U(u,\, v)=\{(p,\, q)\in S\times S:\, d(p,\, q)\leq 3\}=\\
=\Delta_{2}(S)\cup \{(p,\, p+3):\, p=1,\, 2,\, 3\}\cup \{(p+3,\, p):\, p=1,\, 2,\, 3\}\stackrel{def} = \Delta_{3}(S);
\end{multline*}

for every pair of $u$, $v$ such that $4<u\leq 5, 0<v\leq 1$:
\begin{multline*}
U(u,\, v)=\{(p,\, q)\in S\times S:\, d(p,\, q)\leq 4\} = \\
=\Delta_{3}(S)\cup \{(p,\, p+4):\, p=1,\, 2\}\cup \{(p+4,\, p):\, p=1,\, 2\}\stackrel{def} = \Delta_{4}(S);
\end{multline*}

for every pair of $u$, $v$ such that $u>5$, $0<v\leq 1$:
\begin{multline*}
U(u,\, v)=\{(p,\, q)\in S\times S:\, d(p,\, q)\leq 5\}=\\
=\Delta_{4}(S)\cup \{(p,\, p+5):\, p=1\}\cup\{(p+5,\, p):\, p=1\} = S\times S.
\end{multline*}

Also, it is easy to see that $U(u,\, v)=S\times S$ for every $u$, $v$ such that $u>0$, $v>1$.}
\end{Ex}

For any set $Z$ of ordered pairs of positive numbers, i.~e. $Z\subset (0,+\infty)^{2}$, let
\begin{gather*}
\mathcal{N}(Z)=\{N_{p}(u,\, v):\, (u,\, v)\in Z,\, \, p \in S\}\qquad \text{and}\qquad \mathcal{U}(Z)=\{U(u,\, v):\, (u,\,v)\in Z\}.
\end{gather*}

A non-null collection $\{N_{p}\}$ of subsets $N_{p}$ of a set $S$ associated with a point $p\in S$ is a {\itshape family of neighborhoods} for $p$ if each $N_{p}$ contains $p$. Let the family of neighborhoods be associated with each point $p$ of a set $S$.
\begin{itemize}
\item[$V$]\label{typeV}
In this case the set $S$ and the collection of neighborhoods is called the $g$-topological space of type $V$ \cite{Thorp1962}.
\end{itemize}

The {\itshape closure} of a subset $E$ of $S$, written $\overline{E}$, is the set of points $p$ such that each neighborhood of $p$ intersects $E$. The {\itshape interior} of $E$ is the complement of the closure of the complement of $E$. A $g$-topological space $S$ is {\itshape symmetric} if, for every pair of points $p$ and $q$, $p$ is in $\overline{\{q\}}$ iff $q$ is in $\overline{\{p\}}$.

E.~Thorp introduced the following $g$-topologies in a statistical metric space $(S,\, F)$.
\begin{itemize}
\item[$N_{0}$] is type $V$.
\item[$N_{1}$]. For each point $p$ and each neighborhood $U_{p}$ of $p$, there is a neighborhood $W_{p}$ of $p$ such that for each point $q$ of $W_{p}$, there is a neighborhood $U_{q}$ of $q$ contained in $U_{p}$.

\item[$N_{2}$]. For each point $p$ and each pair of neighborhoods $U_{p}$ and $W_{p}$ of $p$, there is a neighborhood of $p$ contained in the intersection of $U_{p}$ and $W_{p}$.
\end{itemize}

The following are various $g$-topologies in a statistical metric space $(S,F)$ defined by E.~Thorp.
\begin{itemize}
\item[$V_{D}$]. If the conditions $N_{0}$ and $N_{2}$ are satisfied, then the collection of neighborhoods on $S$ is called the {\itshape $g$-topology} of type $V_{D}$.

\item[$V_{\alpha}$]. The collection of neighborhoods on $S$ is called the {\itshape $g$-topology} of type $V_{\alpha}$ if the conditions $N_{0}$ and $N_{1}$ are satisfied.

\item[$Top$]. A $g$-topology is a $topology$ if the conditions $N_{0}$, $N_{1}$ and $N_{2}$ are satisfied.

\end{itemize}

Let $S$ be a set and $(P,\, <)$ be a partially ordered set with least element $0$. A {\itshape generalized \'{e}cart} ({\itshape $g$-\'{e}cart} for short) is a mapping
\begin{gather*}
G\colon\, S\times S \to P.
\end{gather*}
If a $g$-\'{e}cart $G$ satisfies $G(p,\, p)=0$ and the set $S$ consists of more than one point, the {\itshape $g$-\'{e}cart $g$-topology} for $S$ is the $g$-topology determined from $G$, and its partially ordered range set $P$, as follows:

For each $f>0$ in $P$ and each $p\in S$, the $f$-sphere for $p$ is a set of the form
\begin{gather*}
N_{p}(f)=\{q\in S:\, G(p,\, q)<f\}.
\end{gather*}
Then for each $p\in S$, the collection of $f$-spheres
\begin{gather*}
\mathcal{N}_{p}(P)=\{N_{p}(f):\, f>0,\, p\in P\}
\end{gather*}
is a family of neighborhoods for $p$.

The $g$-\'{e}cart associated with a statistical metric space $(S,F)$ is the mapping $G$ defined by $G(p,\, q)=G_{pq}$.

\begin{Ex}\label{Npf}{\rm
Let $S=\mathbb{N}$ and $P=\mathbb{N}\cup\{0\}$ be a partially ordered set with the relation $<$ where $\mathbb{N}$ denote the set of all positive integers. Let $A=\{1,\, 2,\, 3\}$ be a subset of $S$. Define
\begin{gather*}
G(p,\, q)=\begin{cases}
1, &\text{if}\,\, p\notin A,\, q\in S,\\
1, &\text{if}\,\, \in A,\, q\notin S,\\
0, &\text{if}\,\, \notin A,\, q\notin S.
\end{cases}
\end{gather*}
and for $p\in A$, $q\in A$ define $G(p,\, q)$ as follows:
\begin{gather*}
G(1,\, 1)=0,\qquad G(1,\, 2)=2,\qquad G(1,\, 3)=3,\\
G(2,\, 1)=4,\qquad G(2,\, 2)=0,\qquad G(2,\, 3)=6,\\
G(3,\, 1)=1,\qquad G(3,\, 2)=2,\qquad G(3,\, 3)=0.
\end{gather*}

$f$-sphere for each $p\in S$ has the following form (A.~A.~Zaitov):

\textbf{Case} $p=1$:
\begin{gather*}
N_{1}(f)=\begin{cases}
\emptyset, &\text{if}\,\, f=0,\\
\{1\}, &\text{if}\,\, f=1,\\
S\setminus \{2,\, 3\}, &\text{if}\,\, f=2,\\
S\setminus \{3\}, &\text{if}\,\, f=3,\\
S, &\text{if}\,\, f\geq 4.
\end{cases}
\end{gather*}

\textbf{Case} $p=2$:
\begin{gather*}
N_{2}(f)=\begin{cases}
\emptyset, &\text{if}\,\, f=0,\\
\{2\}, &\text{if}\,\, f=1,\\
S\setminus \{1,\, 3\}, &\text{if}\,\, 2\leq f\leq 4,\\
S\setminus \{3\}, &\text{if}\,\, 5\leq f\leq 6,\\
S, &\text{if}\,\, f\geq 7.
\end{cases}
\end{gather*}

\textbf{Case} $p=3$:
\begin{gather*}
N_{3}(f)=\begin{cases}
\emptyset, &\text{if}\,\, f=0,\\
\{3\}, &\text{if}\,\, f=1,\\
S\setminus \{2\}, &\text{if}\,\, f=2,\\
S, &\text{if}\,\, f\geq 3.
\end{cases}
\end{gather*}

\textbf{Case} of arbitrary $p\in S\setminus A$:
\begin{gather*}
N_{p}(f)=\begin{cases}
\emptyset, &\text{if}\,\, f=0,\\
S\setminus A, &\text{if}\,\, f=1,\\
S, &\text{if}\,\, f\geq 2.
\end{cases}
\end{gather*}}
\end{Ex}

Given a statistical metric space $(S,\, F)$, for each pair of points $p$ and $r$ in $S$, number $u>0$, the $r$-sphere with center $p$,
$N_{p} (r;\, u)$ is defined to be the sphere
\begin{gather*}
N_{p}(r;\, u)=N(G_{pr}(u))=\{q:\, G_{pq} (u)<G_{pr} (u)\}.
\end{gather*}
The $R$-$g$-topology for $(S,\, F)$ is the structure whose family of neighborhoods at each point $p$ is the collection
\begin{gather*}
\mathcal{N}_{p}(r)=\{N_p (r;\,u):\, r \in S,\, u>0\}.
\end{gather*}

\begin{Ex}{\rm
Consider the $SM$ space $(S,\, F)$ where $S=\mathbb{N}$ and the distribution function
\begin{gather*}
F_{pq}(x)=\begin{cases}
\frac{x}{d(p,\, q)}, &\text{if}\,\, 0<x<d(p,\, q), \,\, d(p,\, q)\neq 0,\\
1, &\text{if}\,\, x\geq d(p,\, q).
\end{cases}
\end{gather*}
where $d(p,\, q)=|q-p|$, $p,\, q \in S$.

Fix $p=1$ and $r=2$ from $S$. Let $x=\frac{1}{4}$. Then $G_{pr}(x)=G_{12}(\frac{1}{4})=1-F_{12}(\frac{1}{4})=0.75$. We have (A.~A.~Zaitov):
$F_{11}(\frac{1}{4})=1>0.25$ and $F_{1q}(\frac{1}{4})=\frac{1}{4}(q-1)\leq 0.25$ for every $q\geq 2$. That is why
\begin{multline*}
N_{1}\left(2;\, \frac{1}{4}\right)=N\left(G_{12}\left(\frac{1}{4}\right)\right) =\\
= \left\{q\in S:\, G_{1q} \left(\frac{1}{4}\right)<0.75\right\} = \left\{q\in S:\, F_{1q}(\frac{1}{4})>0.25\right\}=\{1\}.
\end{multline*}

Note that $N(G_{pr})  = \varnothing$ if $p=r$. Really, for every $u>0$ and $p\in S$ one has $F_{pp}(u)=1$, consequently,
$G_{pp}(u)=0$. Since $0\leq G_{pq}(u)\leq 1$ for all $u>0$, $p,\, q\in S$, there exists no
$q\in S$ such that $G_{pq}(u)<G_{pp}(u)$. Hence $N(G_{pp})=\varnothing$, $p\in S$.}
\end{Ex}

\begin{Rem}{\rm In a $SM$ space $(S,F)$, we use the following notations:
\begin{itemize}
\item[$(a)$] Let $\tau$ denote the $g$-topology of type $V$.
\item[$(b)$] Let $\tau_{D}$ denote the $g$-topology of type $V_{D}$.
\item[$(c)$] Let $\tau_{\alpha}$ denote the $g$-topology of type $V_{\alpha}$.
\item[$(d)$] Let $\tau_{e}$ denote the $g$-\'{e}cart $g$-topology.
\item[$(e)$] Let $\tau_{R}$ denote the $R$-$g$-topology.
\item[$(f)$] Each element in $\mathcal{N}(X)$ is called a $\tau$-neighborhood.
\item[$(g)$] Each element in $\mathcal{N}_{p}(P)$ is called a $\tau_{e}$-neighborhood.
\item[$(h)$] Each element in $\mathcal{N}_{p}(r)$ is called a $\tau_{R}$-neighborhood.
\end{itemize}}
\end{Rem}

\section{Main part}

Let $(S^{1},\, F^{1})$ and $(S^{2},\, F^{2})$ be statistical metric spaces, $(p^{1},\, p^{2})$, $(q^{1},\, q^{2})$, $(r^{1},\, r^{2})\in S^{1}\times S^{2}$. Put
\begin{gather}\label{p.ofd.f.}
F^{12}_{(p^{1},\, p^{2})(q^{1},\, q^{2})} (x) = F^{1}_{p^{1}q^{1}}(x)\cdot  F^{2}_{p^{2}q^{2}}(x).
\end{gather}
Obviously, conditions $(SM\mbox{-}I)$ -- $(SM\mbox{-}III)$ are true. Let us show $(SM\mbox{-}IV)$ is also satisfied. Assume that $F^{12}_{(p^{1},\, p^{2})(q^{1},\, q^{2})}(x)=1$ and $F^{12}_{(q^{1},\, q^{2})(r^{1},\, r^{2})}(y)=1$. These equalities mean that the ``distance'' between $(p^{1},\, p^{2})$ and $(q^{1},\, q^{2})$ less then $x$, and the ``distance'' between $(q^{1},\, q^{2})$ and $(r^{1},\, r^{2})$ less than $y$. Then clearly, that the ``distance'' between $(p^{1},\, p^{2})$ and $(r^{1},\, r^{2})$ less than $x+y$. Consequently, the probability that the ``distance'' between $(p^{1},\, p^{2})$ and $(r^{1},\, r^{2})$ is less than $x+y$, is 1, i.~e. $F^{12}_{(p^{1},\, p^{2})(r^{1},\, q^{2})}(x+y)=1$.

So, (\ref{p.ofd.f.}) is defined correctly. Its tail
\begin{gather*}
G^{12}_{(p^{1},\, q^{1})(p^{2},\, q^{2})}(x) = 1 - F^{12}_{(p^{1},\, q^{1})(p^{2},\, q^{2})}(x)\qquad \mbox{for all}\qquad x \in\mathbb{R}.
\end{gather*}

Now we check the Menger inequality. Let $T(a,\, b) = ab$, $a,\, b\in [0,\, 1]$. Clearly, it satisfies conditions $(T\mbox{-}I)$ -- $(T\mbox{-}V)$. Suppose
\begin{gather*}
F^{1}_{p^{1}r^{1}}(x+y)\geq T(F^{1}_{p^{1}q^{1}}(x),\, F^{1}_{q^{1}r^{1}}(y))\quad \mbox{and}\quad F^{2}_{p^{2}r^{2}}(x+y)\geq T(F^{2}_{p^{2}q^{2}}(x),\, F^{2}_{q^{2}r^{2}}(y)).
\end{gather*}
Then
\begin{align*}
F^{12}_{(p^{1},\, p^{2})(r^{1},\, r^{2})}(x+y) &= F^{1}_{p^{1}r^{1}}(x+y)\cdot F^{2}_{p^{2}r^{2}}(x+y)\geq\\
&\geq T(F^{1}_{p^{1}q^{1}}(x),\, F^{1}_{q^{1}r^{1}}(y))\cdot T(F^{2}_{p^{2}q^{2}}(x),\, F^{2}_{q^{2}r^{2}}(y)) =\\
&=F^{1}_{p^{1}q^{1}}(x)\cdot F^{1}_{q^{1}r^{1}}(y)\cdot F^{2}_{p^{2}q^{2}}(x)\cdot F^{2}_{q^{2}r^{2}}(y)=\\
&=F^{12}_{(p^{1},\, p^{2})(q^{1},\, q^{2})} (x)\cdot F^{12}_{(q^{1},\, q^{2})(r^{1},\, r^{2})} (y) = \\
&= T(F^{12}_{(p^{1},\, p^{2})(q^{1},\, q^{2})}(x),\, F^{12}_{(q^{1},\, q^{2})(r^{1},\, r^{2})} (y)),
\end{align*}
i.~e.
\begin{gather*}
F^{12}_{(p^{1},\, p^{2})(r^{1},\, r^{2})}(x+y)\geq T(F^{12}_{(p^{1},\, p^{2})(q^{1},\, q^{2})}(x),\, F^{12}_{(q^{1},\, q^{2})(r^{1},\, r^{2})} (y)).
\end{gather*}

Now we claim the following statements the proofs each of them consists just directly verification.

\begin{Th}{\rm
Let $(S^{1},\, \tau^{1})$ and $(S^{2},\, \tau^{2})$ be $g$-topological spaces of type $V$. Then $\tau^{1}\times \tau^{2}$ is a $g$-topology of type $V$ on $S_{1}\times S_{2}$.}
\end{Th}

\begin{Th}{\rm
Let $(S^{1},\, \tau^{1}_{D})$ and $(S^{2},\, \tau^{2}_{D})$ be $g$-topological spaces of type $V_{D}$. Then $\tau^{1}_{D}\times \tau^{2}_{D}$ is a $g$-topology of type $V_{D}$ on $S_{1}\times S_{2}$.}
\end{Th}

\begin{Th}{\rm
Let $(S^{1},\, \tau^{1}_{\alpha})$ and $(S^{2},\, \tau^{2}_{\alpha})$ be $g$-topological spaces of type $V_{\alpha}$. Then $\tau^{1}_{\alpha}\times \tau^{2}_{\alpha}$ is a $g$-topology of type $V_{\alpha}$ on $S_{1}\times S_{2}$.}
\end{Th}

\begin{Th}{\rm
Let $(S^{1},\, \tau^{1}_{e})$ and $(S^{2},\, \tau^{2}_{e})$ be $g$-\'{e}cart  $g$-topological spaces. Then $\tau^{1}_{e}\times \tau^{2}_{e}$ is  a $g$-\'{e}cart  $g$-topology on $S_{1}\times S_{2}$.}
\end{Th}

\begin{Th}{\rm
Let $(S^{1},\, \tau^{1}_{e})$ and $(S^{2},\, \tau^{2}_{e})$ be $R$-$g$-topological spaces. Then $\tau^{1}_{e}\times \tau^{2}_{e}$ is  a $R$-$g$-topology on $S_{1}\times S_{2}$.}
\end{Th}

\begin{Rem}{\rm
Note that the Cartesian product of topologies must not be a topology, i.~e. usual topology does not close under product. Unlike usual topology, $g$-topologies are closed with respect to product.

}
\end{Rem}


\begin{thebibliography}{00}
\bibitem{Menger1942}
Menger,~K. Statistical metrics, Proc. Nat. Acad. Sci. USA 28, 535 -- 537 (1942).

\bibitem{Wald1943}
Wald,~A. On a Statistical Generalization of Metric Spaces,  Proc. Nat. Acad. Sci. 29, 196 -- 197 (1943).

\bibitem{Thorp1962}
Thorp,~E.~O. Generalized topologies for statistical metric spaces, Fund. Math., Li, 51, 9 -- 21 (1962).


\bibitem{RenVad2019}
Renukadevi,~V., Vadakasi,~S. On Various g-Topology in Statistical Metric Spaces, Universal Journal of Mathematics and
    Applications, 2(3), 107 -- 115 (2019).

\bibitem{Zaitov2019NUUz}
Zaitov,~A.~A.   Statistical metric spaces, Abstract of the International Conference `Modern Problems of Geometry and Topology and their Applications', Tashkent, Uzbekistam, November 21-23, 2019, P. 89 -- 90.


\end{thebibliography}
\end{document}